\documentclass[11pt]{article}

\usepackage[english]{babel}

\usepackage{authblk}

\usepackage{color}
\usepackage{array}
\usepackage{enumerate}
\usepackage{enumitem}
\usepackage{refcount}
\usepackage{algorithm}
\usepackage{algorithmic}
\usepackage{boxedminipage}
\usepackage{float}
\usepackage[utf8x]{inputenc}
\usepackage{authblk}
 
\usepackage{color}
\usepackage{graphicx}
\usepackage{multirow}
\usepackage{tikz-cd}
\usepackage{amssymb}
\usepackage{caption}
\usepackage{amsmath}
\usepackage{amsfonts}
\usepackage{mathtools}

\newtheorem{thm}{Theorem}[section]
\newtheorem{prop}[thm]{Proposition}
\newtheorem{lem}[thm]{Lemma}

\newtheorem{defn}[thm]{Definition}
\newtheorem{nota}[thm]{Notation}
\newtheorem{exa}[thm]{Example}

\newtheorem{remark}[thm]{Remark}
\newtheorem{conj}[thm]{Conjecture}

\def\proof{{\bf Proof} }

\newcommand{\bC}{\mathcal C} 

\newcommand{\bI}{\mathcal{I}}

\newcommand{\Mod}[1]{\ (\mathrm{mod}\ #1)}

\newcommand{\prfend}{\hbox to7pt{\hfil}
\par\vskip-\baselineskip\hbox to\hsize
{\hfil\vbox {\hrule width6pt height6pt}}\vskip\baselineskip}

\def\dd{\medskip \par \noindent}
\long\def\eatit#1{}

\def\Z{\Bbb Z}
\def\N{\Bbb N}
\def\C{\Bbb C}

\def\P{\Bbb P}

\def\X{\Bbb X}
\def\K{\Bbb K}
\def\Y{\Bbb Y}

\font\tengothic=eufm10
\font\sevengothic=eufm7
\newfam\gothicfam
       \textfont\gothicfam=\tengothic
       \scriptfont\gothicfam=\sevengothic

    \font\tenmsb=msbm10              \font\sevenmsb=msbm7
\newfam\msbfam
       \textfont\msbfam=\tenmsb
       \scriptfont\msbfam=\sevenmsb
\def\Bbb#1{{\fam\msbfam #1}}

\usepackage{xypic}

\DeclareMathOperator{\cha}{char}
\DeclareMathOperator{\im}{Im}
\DeclareMathOperator{\codim}{codim}
\DeclarePairedDelimiter\ceil{\lceil}{\rceil}
\DeclarePairedDelimiter\floor{\lfloor}{\rfloor}

\title{Complete intersections on Veronese surfaces}
\date{June 14, 2022}
\author[1]{Stefano Canino\thanks{stefano.canino@polito.it}}
\author[1]{Enrico Carlini \thanks{enrico.carlini@polito.it}}
\affil[1]{DISMA-Department of Mathematics, Politenico di Torino, Corso Duca degli Abruzzi 24, 10129, Torino}

\begin{document}

\maketitle 

\begin{abstract}
In this paper we describe all possible reduced complete intersection sets of points on Veronese surfaces. We formulate a conjecture for the general case of complete intersection subvarieties of any dimension and we prove it in the case of the quadratic Veronese threefold. Our main tool is an effective characterization of all possible Hilbert functions of reduced subvarieties of Veronese surfaces.
\end{abstract}

\section{Introduction}
Complete intersection subavarities are both a classical both a modern topic of study in Algebraic Geometry.
\dd
In \cite{Commentationes} Euler asked when a sets of points in the plane is the intersection of two curves, that is, using the modern terminology, when a set of points in the plane is a complete intersection. In the same period, Cramer asked similar questions so that this type of questions is presently known as the {\it Cramer-Euler problem}. Euler solution in the case of nine points in the plane gave rise to what are now known as the Cayley-Bacharach Theorems, see \cite{CayleyBacharach}.
\dd
Complete intersections and their algebraic counterpart, regular sequences, play a central role in Commutative Algebra and in Algebraic geometry. Consider, for example, the well known Hartshorne conjecture, stated for the first time by Hartshorne in \cite{Amplesub} and still open, which is probably one of the most studied problems regarding complete intersection. More recently, complete intersections have shown to have unexpected applications. For example, in \cite{SSRF} and \cite{Onthes}, the strength and the slice rank of polynomials are studied using complete intersections. See also \cite{Examplesof} for an application in proving the existence of special families of vector bundles on quartic surfaces of $\P^3$. For a a more exhaustive overview on complete intersections we advise to see \cite{CIAcireale}.
\dd
In this paper,  we consider a generalization of the Cramer-Euler problem: characterize the possible complete intersections lying on a Veronese surface $V_{2,d}$, and more generally on a Veronese variety. We recall that, given non-negative integers $n$ and $d$, a Veronese variety of $\P^N$ is  projectively equivalent to $V_{n,d}\subseteq\P^N$, that is the image of $\P^n$ via the usual Veronese map:
$$\begin{array}{cccc}
\nu_{n,d}: &\P^n  & \rightarrow & \P^N\\
& [x_0,\dots,x_n] & \mapsto & [x_0^d,x_0^{d-1}x_1,\dots,x_n^d]
\end{array}$$
where $N={n+d\choose d}-1$. Notice that for $d=1$ the Veronese surface $V_{2,1}$ is the plane $\P^2$, so that our problem in this special case is exactly the Cramer-Euler problem.
In Theorem \ref{cisuperfici} we completely solve the problem showing that for $d>2$ the only reduced complete intersections of $\P^N$ lying on $V_{2,d}$ are finite sets of either one or two points. While, for the Veronese surface $V_{2,2}\subseteq\P^5$, one also has plane conics and their intersections with suitable hypersurfaces. Moreover, in Theorem \ref{cicurve} we show that, except for the case $d=2$, the only complete intersections lying on rational normal curves $V_{1,d}$ are the trivial ones, that is one single point or the set of two points. The case $V_{1,2}$, that is of a plane conic, is different. In fact, by cutting with any properly chosen curve, one will produce a complete intersection set of points. 
\dd
Inspired by these evidences we formulate Conjecture \ref{theconjecture}: the only reduced complete intersections of $V_{n,d}$, $d\geq 3$, are finite sets of either one or two points while for $d=2$ one also has plane conics and their intersections with suitable hypersurfaces. We also checked the validity of the conjecture for $V_{3,2}$, see Proposition \ref{propositionV32}.
\dd
In order to prove our main result, Theorem \ref{cisuperfici}, we characterize the possible Hilbert functions of reduced subvarieties of Veronese varieties in in Theorem \ref{strumento}. In other words, we characterize all possible Hilbert functions of radical ideals in the Veronese ring $\C[x_0,\dots,x_N]/\bI(V_{n,d})$, where $\bI(V_{n,d})$ is the defining ideal of $V_{n,d}$. Beyond their application to the proof of our theorem, Hilbert functions play a central role in Commutative Algebra and in Algebraic Geometry, for example see \cite{StanleyHF},  \cite{ThegeometryofHilbertfunctions}, and \cite{BigattiGeramitaMigliore}. Indeed, in recent times Hilbert functions have also been used as tools in other fields. For example, in the study of Waring rank, that is the tensor rank for symmetric tensors, in the paper \cite{ThesolutionofWproblem}. Another example is the study of Strassen's Conjecture, a crucial conjecture in complexity theory, now proved to be false in general \cite{Counterexamplesstrassen}, but still open in the relevant case of symmetric tensors, see \cite{Symmetrictensors} and \cite{OnComons}. As a last example, we also mention the study of the identifiability of tensors, which plays a crucial role in Algebraic Statistic, see \cite{Realidentifiability},\cite{Ontheidentifiability} and \cite{OnthenumberofWaring}.
\dd 
In this paper, se first characterize the Hilbert functions of  reduced subvarieties of  $V_{n,d}$. Thus, we generalise 0-sequences and differentiable 0-sequences introduced in \cite{GMR}, see Definition \ref{zeroseq} and Definition \ref{diffzeroseq}. Successively, we give more effective characterization for the case of the rational normal curves $V_{1,d}$ in Theorem \ref{hfcurve}, so recovering a classical result, and for the case of the surfaces $V_{2,d}$ in Theorem \ref{main}. In \cite[Theorem 4.5]{Hilbertschemes} a similar characterization is given in the case of subschemes of $V_{n,d}$, that is in the case of any ideal in a Veronese ring. However, we consider these characterization not to be enough effective for our purposes. In fact, given a candidate Hilbert function, one has to solve a kind of interpolation problem to decide whether the given function is the Hilbert function of a reduced subvariety, or subscheme, of $V_{n,d}$, see Remark \ref{interpolazione}.
\dd
More precisely, the paper is structured as follows.
In \S \ref{SectionFirstResult}, we recall some basic notions needed in the paper and we study the relationship between $\bI(\nu_{n,d}^{-1}(\X))$ and $\bI(\X)$, where $\X\subseteq V_{n,d}$ is a reduced subvariety of $V_{n,d}$. Using this we characterize the Hilbert functions of subvarieties lying on $V_{n,d}$ and introduce $d$-sequences and differentiable $d$-sequences. Moreover, we use these results to study the Hilbert functions of divisors of $V_{n,d}$ and the Hilbert functions of finite sets of reduced points lying on a rational normal curve $V_{1,d}$, recovering the classical results.
In \S \ref{The case of surfaces}, we apply Theorem \ref{Cond} to the special case of Veronese surfaces getting a more explicit characterization for the Hilbert functions of reduced subvarieties of Veronese surfaces in  Theorem \ref{main}.
In \S \ref{An application to complete intersections}, we use the previously introduced tools to prove our main theorem characterizing all possible complete intersections lying on Veronese surfaces.
In \S \ref{More results and open problems}, we study the reduced complete intersection lying on $V_{3,2}$ and state Conjecture \ref{theconjecture} about reduced complete intersections lying on Veronese varieties.
\dd 
{\bf Acknowledgments}: the authors wish to thank Ciro Ciliberto, Luca Chiantini, Alessandro \linebreak Gimigliano, and Giorgio Ottaviani for some useful discussions. Moreover, the authors acknowledges the kindness of Marvi Catalisano and Juan Migliore in parting from a similar project which was never completed. The authors are members of GNSAGA of INDAM and were was supported by MIUR grant Dipartimenti
di Eccellenza 2018-2022 (E11G18000350001).

\section{Basic notions and first results}\label{SectionFirstResult}
In this section we introduce the needed basic notions and some preliminary results including a complete characterization of the possible Hilbert functions of reduced subvarieties of a Veronese variety.
Let $\K$ be an algebraic closed field with $\cha{\K}=0$ and  $R=\K[x_0,\dots,x_n]$ the coordinate ring of $\P^n$. If $d$ is a positive integer, we set $N={n+d\choose d}-1$ and we denote the coordinate ring of $\P^N$ by $S=\K[y_0,\dots,y_N]$.
\dd 
If $M$ is a graded module, we denote by $M_p$ the $p$-th graded component of $M$, so that $M=\bigoplus_{p=0}^\infty M_p$. In particular, since an ideal $I\subseteq R$ is an $R$-module, we use the same notation $I_p$ for the $p$-th graded component of $I$.
\dd 
For a graded module $M$, we indicate its Hilbert function by $H_M(t)=\dim M_t$ and its first difference function by $\Delta H_M(t)=H_M(t)-H_M(t-1)$, setting $H_M(-1)=0$. If $\X\subseteq\P^n$ is a reduced projective scheme we denote by $\bI(\X)$ its associated ideal and by $H_\X(t)$, its Hilbert function, that is the Hilbert function of the module $R/\bI(\X)$.
\begin{nota}\rm
We will  denote by $\alpha_0,\dots,\alpha_N$ the multi-indeces of $\left\lbrace 0,1,\dots,d\right\rbrace ^{n+1}$ such that $|\alpha_j|=d $ ordered by the usual lexicographic order. Moreover, if $\alpha_i=(\alpha_{i0},\dots,\alpha_{in})$, we set 
$$ \qquad\underline{x}^{\alpha_i}=x_0^{\alpha_{i0}}\dots x_n^{\alpha_{in}}.$$
\end{nota}
\begin{defn}\rm 
For each $n,d\in\N_{>0}$ we define the \textit{$(n,d)$-Veronese embedding}
$$\begin{array}{cccc}
\nu_{n,d}: &\P^n  & \rightarrow & \P^N\\
& [x_0,\dots,x_n] & \mapsto & [\underline{x}^{\alpha_0},\dots,\underline{x}^{\alpha_N}]
\end{array}$$
 and   the $(n,d)$\textit{-Veronese variety} $V_{n,d}:=\nu_{n,d}(\P^n)$.
\end{defn}
The following lemma will be very useful.
\begin{lem}\label{Imphi}\rm 
Let us consider the graded morphism
$$\begin{array}{cccc}
\varphi_d: & S & \rightarrow & R\\
& a & \mapsto & a\\
& y_i & \mapsto & \underline{x}^{\alpha_i}
\end{array}$$
for all $a\in\K$, and for  $i\in\{0,\dots,N\}$. Then the following hold:
\begin{enumerate}
	\item $\ker\varphi_d=\bI(V_{n,d})$.
	\item $\im\varphi_d=\bigoplus_{\ell=0}^\infty R_{\ell d}$ and, in particular, $\varphi_d(S_t)=R_{td}$.
	\item If $\X\subseteq V_{n,d}$, then $(\bI(\nu_{n,d}^{-1}(\X)))_{td}=\varphi_d(\bI(\X)_t)$. In particular
$$\varphi_d(\bI(\X))=\bigoplus_{s=0}^\infty (\bI(\nu_{n,d}^{-1}(\X)))_{sd}.$$
\item If $\X$ is a subvariety of $V_{n,d}$ and we set $\Y=\nu_{n,d}^{-1}(\X)$, then 
$$H_{\X}(t)=H_{\Y}(td)\;\forall\;t\geq0.$$
\end{enumerate}
\end{lem}
\proof The proof of (1) is trivial. To prove (2) it is enough to note that any monomial of degree $td$ can be written as product of monomials of degree $d$. The proof of (3) is a straightforward check of a double inclusion. The proof of (4) follows from the chain of graded isomorphisms
$$S/\bI(\X)\cong{S/\bI(V_{n,d})\over \bI(\X)/\bI(V_{n,d})}\cong {\psi_d(S/\bI(V_{n,d}))\over \psi_d(\bI(\X)/\bI(V_{n,d}))}={ \bigoplus_{s=0}^\infty R_{sd}\over \varphi_d(\bI(\X))}$$
where $\psi_d:S/\bI(V_{n,d})  \rightarrow  {\bigoplus_{s=0}^\infty R_{sd}}$ is the canonical isomorphism induced by $\varphi_d$.
\prfend
\begin{remark}\label{Boh}\rm
\noindent Since $\varphi_d$ is a ring homomorphism, $\varphi_d(\bI(\X))$ is an ideal of $\im\varphi_d=\bigoplus_{s=0}^\infty R_{sd}$, but it is not an ideal of $R$. Nevertheless, one has $$(\varphi_d(\bI(\X))R)_{td}=(\varphi_d(\bI(\X)))_{td}.$$
\end{remark}
\begin{remark}\label{hfveronese}\rm 
In the notations of Lemma \ref{Imphi}, if we choose $\X=V_{n,d}$ then we have $\Y=\P^n$ and thus we get that the Hilbert Function of $V_{n,d}$ is $H_{V_{n,d}}(t)={n+td\choose n}$.
\end{remark}
The following theorem is an immediate consequence of Lemma \ref{Imphi}.
\begin{thm}\label{Cond}
Let $h(t):\N\to\N$ be the Hilbert function of a projective variety in $\P^N$. Then there exists $\X\subseteq V_{n,d}\subseteq\P^N$ such that $H_\X(t)=h(t)$ if and only there exists $k(t):\N\to\N$ Hilbert function of a projective variety in $\P^n$ such that $h(t)=k(dt)$.
\end{thm}
\begin{remark}\rm
We note that, if a subvariety $\Y$ has Hilbert function satisfying the conditions of Theorem \ref{Cond},  this does not mean that $\Y$ lies on a Veronese variety. The theorem only guarantees that there exists {\em some} subvariety $\X$ of a Veronese variety having the same Hilbert function of $\Y$. Consider, for example, seven generic points in $\P^3$. By genericity they do not lie on a $V_{1,d}$, that is they do not lie on rational normal curve, but their Hilbert function satisfies the hypothesis of the theorem.
\end{remark}
In the case of divisors of $V_{n,d}$ we can be more explicit.
\begin{prop}\label{hfdivisors}
If  $\X$ is a divisor of $V_{n,d}$, let $\deg\X=de$, then
$$H_\X(t)=\begin{cases}
{n+dt\choose n}, &\text{if } t\leq\floor[\big]{{e-1\over d}} \\
{n+dt\choose n}-{n+dt-e\choose n}, &\text{if } t\geq \floor[\big]{{e-1\over d}}+1
\end{cases}$$
\end{prop}
\proof
Since $\X$ is a divisor and $\deg\X=de$, there exists a (unique) hypersurface $\Y:F=0$ of degree $e$ in $\P^n$ such that $\nu_{n,d}(\Y)=\X$. For each $t\in\N$ we have the following short exact sequence
$$\begin{tikzcd}
0 \arrow[r] & R(-e)_t \arrow[r, "F"] & R_t \arrow[r,"\pi"] & (R/(F))_t \arrow[r,] & 0
\end{tikzcd}$$
and as a consequence, we get $H_{\Y}(t)=\dim (R/(F))_t=\dim R_t-\dim R(-e)_t$ and this ends the proof.
\begin{remark}\label{hfcurve}\rm 
As special case of divisor we can consider $\X\subseteq V_{1,d}\subseteq \P^d$ a finite set of $s$ reduced points on the rational normal curve of degree $d$. Using Proposition \ref{hfdivisors} we get the well known result:$$H_\X(t)=\begin{cases}
	dt+1, &\text{if } t\leq\floor[\big]{{s-2\over d}}\\
	s, &\text{if } t\geq\floor[\big]{{s-2\over d}}+1
	\end{cases}.$$
\end{remark} 
We conclude recalling the following definitions from \cite{GMR}.
\begin{defn}\label{zeroseq}\rm
A sequence of non-negative integers $(c_t)_{t\in\N}$ is called a {\it $0$-sequence} if $c_0=1$ and $c_{t+1} \leq c^{\langle t\rangle}$ for all $t\geq 1$.
\end{defn}
\begin{defn}\label{diffzeroseq}\rm
Let $(b_t)_{t\in\N}$ be a 0-sequence. Then $(b_t)_{t\in\N}$ is {\it differentiable} if the difference sequence $(c_t)_{t\in\N},c_t=b_t-b_{t-1}$ is again a $0$-sequence (where $b_{-1}=0$).
\end{defn}
In \cite{GMR} the authors completely characterized the Hilbert functions of reduced varieties with the following theorem.
\dd
\begin{thm}[\cite{GMR},Theorem 3.3]\label{GMRosso}
Let $\K$ be an infinite field and let $(b_t)_{t\in\N}$ be a differentiable 0-sequence with $b_1=n+1$. There is a radical ideal $I$ in $\K[x_0,\dots,x_n]$ such that $(b_t)_{t\in\N}$ is the Hilbert function of $\K[x_0,\dots,x_n]/I$.
\end{thm}
 The complete characterization of Theorem \ref{GMRosso} can be rephrased: if $S$ is a sequence of non-negative integers, then there exists a reduced $\K$-algebra $A$ such that $S$ is the Hilbert function of $A$ if and only if $S$ is a differentiable $0$-sequence.
\dd
Theorem \ref{Cond} suggests us to extend Definition \ref{zeroseq} and \ref{diffzeroseq} as follows.

\begin{defn}\label{dseq}\rm
A $0$-sequence  $(b_t)_{t\in\N}$ is called {\it d-sequence} if there exists a $0$-sequence $(c_t)_{t\in\N}$ such that $b_t=c_{(d+1)t}$.
\end{defn}

\begin{defn}\label{diffdseq}\rm
A 0-sequence $(b_t)_{t\in\N}$ is called {\it differentiable d-sequence} if there exists a differentiable $0$-sequence $(c_t)_{t\in\N}$ such that $b_t=c_{(d+1)t}$.
\end{defn}

\begin{remark}\rm
We note that that a differentiable $d$-sequence is necessarily a differentiable $0$-sequence.
\end{remark}
We can now rephrase Theorem \ref{Cond} as follows:
\begin{thm}\label{strumento}
	Let $(h_t)_{t\in\N}$ be a sequence of non-negative integers such that $h_0=1$ and $h_1=N+1$. There exists a projective variety $\X\subseteq V_{n,d}\subseteq \P^N$ such that $H_{\X}(t)=h_t$ if and only if $(h_t)_{t\in\N}$ is a differentiable $(d-1)$-sequence.
\end{thm}
It is natural to ask for an effective characterization of $d$-sequences similar to the one of Theorem \ref{GMRosso}. The question does not have an answer in general yet, nevertheless one can give an answer in a special case using our results. In the case $(d-1)$-sequences with $h_1={d+2\choose2}$, such a characterization can be easily produced using Theorem \ref{main} in the next section.

\section{Hilbert functions of points on Veronese surfaces}\label{The case of surfaces}

In this section we focus our attention on the case of Veronese surfaces $V_{2,d}$. In particular, we give an effective characterization of the Hilbert function of any reduced subvariety of $V_{2,d}$ in Theorem \nolinebreak\ref{main}.
\begin{nota}\rm
Given $d,t,s\in\N$ such that $s\geq d^2t+{d(d+3)\over 2}$ we define the following two functions:
$$\mu_1(d,t,s):=d^2t+{d(d+3)\over 2}-s$$
$$ \mu_2(d,t,s):=\begin{cases}
\floor[\big]{{2d(t+1)+3-\sqrt{1+8\mu_1(d,t,s)}\over 2}}, &\text{if } 1\leq \mu_1(d,t,s)\leq {d+1\choose 2}\\
dt-n, & \text{if } {d+1\choose 2}+dn<\mu_1(d,t,s)\leq {d+1\choose 2}+d(n+1), 0\leq n\leq dt
\end{cases}$$
\end{nota}
We begin with a technical result.
\begin{lem}\label{tec}
Let $d,t\in\N$. Consider a function $h:\{1,2,\dots,d\}\to\N$ such that there exists \linebreak $i_0\in\{1,2,\dots,d\}$ with the properties
\begin{enumerate}
	\item $h(i)=dt+i+1$ for each $1\leq i\leq i_0-1$;
	\item $h(i)\geq h(i+1)$ for each $i_0\leq i\leq d-1$.
\end{enumerate}
Then $h(d)\leq\mu_2(d,t,\sum_{i=1}^dh(i))$ and moreover the inequality is sharp.
\end{lem}
\proof
See Appendix A.
\prfend
\noindent In the following proposition we characterise Hilbert functions of reduced points in $\P^{{d(d+3)\over 2}}$ which arise from Hilbert functions of reduced points in $\P^2$ by sampling with steps of length $d$.
\begin{prop}\label{Tappabuchi}
		Let us consider a finite set of reduced points $\X\subseteq\P^{{d(d+3)\over 2}}$ and set $$t_1=\max\left\lbrace t\;|\; H_{\X}(t)=H_{V_{2,d}}(t)\right\rbrace,\quad t_2=\min\left\lbrace t\;|\; H_{\X}(t)=|\X|\right\rbrace. $$
		Then $H_{\X}(t)$ is a $(d-1)$-sequence if and only if the following conditions hold 
		\begin{enumerate}[label=\roman*.]
			\item\label{conduno}$$
			\mu_2(d,t_1,\Delta H_{\X}(t_1+1))\geq \ceil[\bigg]{{\Delta H_{\X}(t_1+2)\over d}};$$
			\item\label{conddue} For all $t_1+2\leq t\leq t_2-1$ $$
			\floor[\bigg]{{\Delta H_{\X}(t)\over d}}\geq \ceil[\bigg]{{\Delta H_{\X}(t+1)\over d}}.$$
		\end{enumerate}
	\end{prop}
\proof
 First, we assume that there exists $\Y\subseteq\P^2$ such that $H_{\Y}(dt)=H_{\X}(t)$. Let us set 
$$H_{\X}(t)=h_t, \quad \Delta H_{\X}(t)=\Delta h_t,$$
and
$$H_{\Y}(t)=k_t,\quad \Delta H_{\Y}(t)=\Delta k_t.$$
Since $(k_t)_{t\in\N}$ is the Hilbert function of a finite set of reduced points of $\P^2$, then (see \cite{Extremal}, Proposition 1.1) there exists $t'\in\N$ such that the following conditions hold
\begin{itemize}
		\item $k_t={2+t\choose t}$ for all $t\leq t'$ and $k_t<{2+t\choose t}$ for all $t>t'$, and thus $\Delta k_t=t+1$ for all $t\leq t'$;
	\item $\Delta k_t\geq\Delta k_{t+1}$ for all $t>t'$, and $\Delta k_t$ is eventually equal to $0$.
\end{itemize}
Note that 
$$\Delta h_t=h_t-h_{t-1}=k_{dt}-k_{d(t-1)}=k_{dt}+\sum_{i=1}^{d-1}(k_{dt-i}-k_{dt-i})-k_{dt-d}=\sum_{i=0}^{d-1}(k_{dt-i}-k_{dt-(i+1)})=\sum_{i=0}^{d-1}\Delta k_{dt-i}$$
and thus
\begin{equation}\label{primdiff}
\Delta h_{t+1}=\sum_{i=0}^{d-1}\Delta k_{d(t+1)-i}=\sum_{i=1}^d\Delta k_{dt+i}\;.
\end{equation}
Thus,  Remark \ref{hfveronese} yields that
$$\sum_{i=1}^{d}\Delta k_{dt+i}=\Delta h_{t+1}={2+d(t+1)\choose 2}-{2+dt\choose 2}=d^2t+{d(d+3)\over 2}\text{ for all }t\leq t_1-1\;.$$
Since $\Delta k_t\leq t+1$ for all $t$ and
$$\sum_{i=1}^d (dt+i+1)=d^2t+{d(d+3)\over 2}$$
it follows that
$$\sum_{i=1}^{d}\Delta k_{dt+i}=d^2t+{d(d+3)\over 2}\text{ for all } t\leq t_1-1$$
if and only if
$$ \Delta k_t=t+1\text{ for all } t\leq dt_1$$
and thus $t'\geq dt_1$.  Moreover, since $h_t<{2+dt\choose 2}$ for $t>t_1$, the same argument shows that  
$$\sum_{i=1}^{d}\Delta k_{dt+i}<d^2t+{d(d+3)\over 2}\text{ for all } t\geq t_1$$
and for $t=t_1$ we have
$$\sum_{i=1}^{d}\Delta k_{dt_1+i}<d^2t_1+{d(d+3)\over 2}\;.$$
As a consequence, there exists a minimum $i_0\in\{1,2,\dots,d\}$ such that $\Delta k_{dt_1+i_0}<dt_1+i_0+1$ and therefore $t'\leq dt_1+d-1=d(t_1+1)-1$.
It follows that
\begin{equation}\label{A}
\Delta k_{dt_1+i}=dt_1+i+1
\end{equation}
for $1\leq i\leq i_0-1$
and
$$\Delta k_{dt_1+i_0}\geq\Delta k_{dt_1+i_0+1}\geq\dots\geq\Delta k_{dt_1+i_0+a}=0$$
for some $a\in\N$. Moreover
\begin{equation}\label{B}
\Delta k_{dt_1+i_0}\geq\Delta k_{dt_1+i_0+1}\geq\dots\geq\Delta k_{dt_1+d}
\end{equation}
and
\begin{equation}\label{C}
\Delta k_{dt+1}\geq\Delta k_{dt+2}\geq\dots\geq \underbrace{\Delta k_{dt+d+1}}_{=\Delta k_{d(t+1)+1}}
\end{equation}
for each $t\geq t_1+1$. Now, for a fixed $t\geq t_1+1$, by \eqref{C} we have
\begin{equation*}
(d-1)\Delta k_{dt+1}\geq\sum_{i=2}^d\Delta k_{dt+i},
\end{equation*}
and
\begin{equation*}
(d-1)\Delta k_{dt+d}\leq\sum_{i=1}^{d-1}\Delta k_{dt+i}.
\end{equation*}
Using formula \eqref{primdiff}, we obtain 
$$ \Delta k_{dt+1}=\Delta h_{t+1}-\sum_{i=2}^d\Delta k_{dt+i}\geq\Delta h_{t+1}-(d-1)\Delta k_{dt+1}$$
and hence
$$ \Delta k_{dt+1}\geq {\Delta h_{t+1}\over d}$$
and similarly
$$\Delta k_{dt+d}=\Delta h_{t+1}-\sum_{i=1}^{d-1}\Delta k_{dt+i}\leq\Delta h_{t+1}-(d-1)\Delta k_{dt+d}$$
and hence
$$\Delta k_{dt+d}\leq {\Delta h_{t+1}\over d} .$$
Moreover, since $(h_t)_{t\in\N}$ and $(k_t)_{t\in\N}$ are integer valued, we have that
\begin{gather*}
 \Delta k_{dt+1}\geq\ceil[\bigg] {{\Delta h_{t+1}\over d}},\qquad  \Delta k_{dt+d}\leq\floor[\bigg]{{\Delta h_{t+1}\over d}}
\end{gather*}
for each $t\geq t_1+1$. By \eqref{A} and \eqref{B} it follows that the function
$$\begin{array}{cccc}
\Delta k_{dt_1+i}: & \left\lbrace 1,2,\dots,d\right\rbrace & \to & \N \\
& i & \mapsto & \Delta k_{dt_1+i}\end{array}$$
satisfies the hypothesis of Lemma \ref{tec}. Hence, we get that
$$\mu_2(d,t_1\Delta h_{t_1+1})\geq\Delta k_{dt_1+d}\geq\Delta k_{d(t_1+1)+1}\geq\ceil[\bigg] {{\Delta h_{t_1+2}\over d}} $$
thus proving condition \eqref{conduno}. Finally, using \eqref{C} we get that
$$
\Delta k_{dt+d}\geq \Delta k_{dt+d+1}=\Delta k_{d(t+1)+1}$$
for each $t\geq t_1+1$, and hence 
$$\floor[\bigg]{{\Delta h_{t+1}\over d}}\geq \ceil[\bigg]{{\Delta h_{t+2}\over d}}$$
for each $t\geq t_1+1$. We notice that this inequality is always verified for $t\geq t_2$ since
$$\ceil[\bigg]{{\Delta h_{t+1}\over d}}=0$$
for each $t\geq t_2$. Hence also  condition \eqref{conddue} is proved.
\dd 
Now we assume that conditions \eqref{conduno} and \eqref{conddue} hold and we prove that there exists $\Y\subseteq\P^2$ such that $H_{\Y}(dt)=h_t$. Since 
$$k_t=\sum_{i=0}^t\Delta k_t,$$
we can construct the Hilbert function $(k_t)_{t\in\N}$ by its first difference $\Delta k_t$. For each $t$ let $e_t$ be the unique integer such that 
$$\Delta h_t\equiv e_t\Mod{d},\quad 0\leq e_t\leq d-1.$$
We define $\Delta k_t$ as follows:
\begin{itemize}
	\item if $0\leq t\leq dt_1$ we set $\Delta k_{t}=t+1$;
	\item if $dt_1+1\leq t\leq d(t_1+1)$ we construct $\Delta k_t$ according to Lemma \ref{tec};
	\item if $t\geq d(t_1+1)+1$ we set
	\begin{gather*}
	\Delta k_{dt+i}=\ceil[\bigg]{{\Delta h_{t+1}\over d}},\quad t\geq t_1+1,\;1\leq i\leq e_{t+1}\\
	\Delta k_{dt+i}=\floor[\bigg]{{\Delta h_{t+1}\over d}},\quad t\geq t_1+1,\;e_{t+1}+1\leq i\leq d.
	\end{gather*}
\end{itemize}
Under our assumptions, this choice guarantees that $\Delta k_t$ is the first difference function of a set of reduced points in $\P^2$ (see \cite{Extremal} Proposition 1.1). Moreover, for $t\leq t_1-1$, we have
$$\sum_{i=1}^d\Delta k_{dt+i}=\sum_{i=1}^d (dt+i+1)=d^2t+{d(d+3)\over2}=\Delta h_{t+1},$$
and, for $t=t_1$ we have by construction $$\sum_{i=1}^d\Delta k_{dt_1+1}=\Delta h_{t_1+1}$$
while for $t\geq t_1+1$ we have
\begin{align*}
\sum_{i=1}^d \Delta k_{dt+i} & =\sum_{i=1}^{e_{t+1}}\Delta k_{dt+i}+\sum_{i=e_{t+1}+1}^d\Delta k_{dt+i}=e_{t+1}\ceil[\bigg]{{\Delta h_{t+1}\over d}}+(d-e_{t+1})\floor[\bigg]{{\Delta h_{t+1}\over d}}=  \\
& = e_{t+1}\left( {\Delta h_{t+1}-e_{t+1}\over d}+1\right)+(d-e_{t+1})\left( {\Delta h_{t+1}-e_{t+1}\over d}\right)=\Delta h_{t+1}.  
\end{align*}
Hence we have
$$k_{dt}=\sum_{i=0}^{dt}\Delta k_t=\Delta k_0+\sum_{i=0}^{t-1}\left( \sum_{j=1}^d\Delta k_{di+j}\right)=\Delta h_0+\sum_{i=0}^{t-1}\Delta h_{i+1}=\sum_{i=0}^t\Delta h_t=h_t $$
for all $t$. This concludes the proof.
\prfend

\noindent We now give our effective characterization of the Hilbert functions of reduced sets of points on Veronese surfaces.

\begin{thm}\label{main}
Let $(h_t)_{t\in\N}$ be the Hilbert function of a finite set of $m$ reduced points in $\P^{{d(d+3)\over 2}}$ and set
$$t_1=\max\left\lbrace t\;|\; h(t)=H_{V_{2,d}}(t)\right\rbrace\qquad t_2=\min\left\lbrace t\;|\; h(t)=m\right\rbrace. $$ 
Then there exists $\X\subseteq V_{2,d}\subseteq \P^N$, $|\X|=m$ such that $H_{\X}(t)=h_t$ if and only if the following conditions hold 
\begin{enumerate}[label=\roman*.]
	\item$$
	\mu_2(d,t_1,\Delta h_{t_1+1})\geq \ceil[\bigg]{{\Delta h_{t_1+2}\over d}};$$
	\item For all $t_1+2\leq t\leq  t_2-1$ $$
	\floor[\bigg]{{\Delta h_t\over d}}\geq \ceil[\bigg]{{\Delta h_{t+1}\over d}}.$$
\end{enumerate}
\end{thm}
\proof
It follows from Theorem \ref{Cond} and Proposition \ref{Tappabuchi}.
\prfend
\noindent Let us see now two explicit instances of use of Theorem \ref{main} in determining whether a given Hilbert function can be realised as the Hilbert function of a subvariety of a Veronese surface.
\begin{exa}\label{incompletabile}\rm
Let us consider the sequence $(h_t)_{t\in\N}$ defined as follows
\begin{table}[H]
\centering
\begin{tabular}{c|cccccccccccc}
$t$    & 0 & 1  & 2   & 3   & 4   & 5   & 6   & 7    & 8    & 9    & 10   & 11 \\ \hline \\
$h_t$  & 1 & 36 & 120 & 253 & 435 & 666 & 946 & 1256 & 1531 & 1744 & 1956 & 2022 
\end{tabular}
\end{table}
\noindent and $h_t=2022$ for $t\geq12$. It is easy to check, using Theorem 3.3. in \cite{GMR}, that this is the Hilbert function of a set of 2022 reduced points in $\P^{35}$. We ask whether there exists $\X\subseteq V_{2,7}\subseteq \P^{35}$ such that $H_{\X}(t)=h_t$ for all $t\geq 0$. To answer we use Theorem \ref{main}. First we determine $t_1$ and $t_2$. Since the Hilbert function of $V_{2,7}$ is $H_{V_{2,7}}(t)={2+7t \choose 2}$, we have that
\begin{table}[H]
\centering
\begin{tabular}{c|ccccccccccccc}
$t$           & 0 & 1  & 2   & 3   & 4   & 5   & 6   & 7    & 8    & 9    & 10   & 11   & 12  \\ \hline \\
$H_{V_{2,7}}$ & 1 & 36 & 120 & 253 & 435 & 666 & 946 & 1275 & 1653 & 2080 & 2556 & 3081 & 3655
\end{tabular}
\end{table}
\noindent so that $t_1=6$ and $t_2=11$.  To determine $\mu_1(7,6,\Delta h_{t_1+1})$ we compute $\Delta h_{t_1+1}$. We have that
\begin{table}[H]
\centering
\begin{tabular}{c|ccccccccccccc}
$t$          & 0 & 1  & 2  & 3   & 4   & 5   & 6   & 7   & 8   & 9   & 10  & 11 & 12 \\ \hline \\
$\Delta h_t$ & 1 & 35 & 84 & 133 & 182 & 231 & 280 & 310 & 275 & 213 & 212 & 66 & 0
\end{tabular}
\end{table}
\noindent and thus $\mu_1(7,6,310)=7^2\cdot 6+{7(7+3)\over 2}-310=19$. Finally, since $19\leq {7+1\choose 2}=28$, we get
$$\mu_2(7,6,310)=\floor[\bigg]{{2\cdot 7(6+1)+3-\sqrt{1+8\cdot 19}\over 2}}=44.$$
To check conditions \eqref{conduno} and \eqref{conddue} of Theorem \ref{main}, we compute $\floor[\big]{{\Delta h_t\over 7}}$ and $\ceil[\big]{{\Delta h_t\over 7}}$ obtaining the following table
\begin{table}[H]
\centering
\begin{tabular}{c|ccccccccccccc}
$t$          & 0 & 1 & 2  & 3  & 4  & 5  & 6  & 7  & 8  & 9  & 10 & 11 & 12  \\ \hline \\
$\ceil[\big]{{\Delta h_t\over 7}}$ & 1 & 5 & 12 & 19 & 26 & 33 & 40 & 45 & 40 & 31 & 31 & 10 & 0  \\ \\
$\floor[\big]{{\Delta h_t\over 7}}$           & 0 & 5 & 12 & 19 & 26 & 33 & 40 & 44 & 39 & 30 & 30 & 9  & 0
\end{tabular}
\end{table}
\noindent Since $\mu_2(7,6,310)=44$ and $\ceil[\big]{{\Delta h_8\over 7}}=40 $ condition \eqref{conduno} is satisfied. However condition \eqref{conddue} is not satisfied for $t=9$ and hence such an $\X$ does not exist.
\end{exa}
\begin{exa}\label{completabile}\rm 
Now, we consider the sequence $(h_t)_{t\in\N}$ defined as follows
\begin{table}[H]
\centering
\begin{tabular}{c|cccccccccccc}
$t$    & 0 & 1  & 2   & 3   & 4   & 5   & 6   & 7    & 8    & 9    & \bf{10}   & 11  \\ \hline \\
$h_t$ & 1 & 36 & 120 & 253 & 435 & 666 & 946 & 1256 & 1531 & 1744 & \bf{1915} & 2022
\end{tabular}
\end{table}
\noindent and $h_t=2022$ for $t\geq 12$; note that this function coincides with the one of the previous example, but for $t=10$. We ask whether there exists $\X\subseteq V_{2,7}\subseteq \P^{35}$ such that $H_{\X}(t)=h_t$ for all $t\geq 0$. As in the previous example we have $t_1=6$ and $t_2=11$. Moreover, we get
\begin{table}[H]
\centering
\begin{tabular}{c|ccccccccccccc}
$t$                                 & 0 & 1  & 2  & 3   & 4   & 5   & 6   & 7   & 8   & 9   & 10  & 11  & 12  \\ \hline \\
$\Delta h_t$                       & 1 & 35 & 84 & 133 & 182 & 231 & 280 & 310 & 275 & 213 & 201 & 107 & 0   \\ \\
$\ceil[\big]{{\Delta h_t\over 7}}$ & 1 & 5  & 12 & 19  & 26  & 33  & 40  & 45  & 40  & 31  & 29  & 16  & 0   \\ \\
$\floor[\big]{{\Delta h_t\over 7}}$  & 0 & 5  & 12 & 19  & 26  & 33  & 40  & 44  & 39  & 30  & 28  & 15  & 0 
\end{tabular}
\end{table}
\noindent and thus $\mu_1(7,6,310)=19$ and $\mu_2(7,6,310)=44$. Thus, condition \eqref{conduno} is satisfied and condition \eqref{conddue} is satisfied for $t=8,9,10$. Hence such an $\X$ exists.
\end{exa}
\begin{remark}\label{interpolazione}\rm
One could deal with both the previous examples  without using Theorem \ref{main}, but just using Macaulay's inequalities for Hilbert functions. However, this requires a trial and error approach. In fact, since $d=7$ one should try to fill step by step six gaps between $h_i$ and $h_{i+1}$ satisfying Macaulay's inequalities. Thus, for each possible choice at each step, one should compute the appropriate binomial expansion. Moreover, one should check also that the first difference function of the constructed Hilbert function is still an Hilbert function: doing this will require a very large number of computations.
\end{remark}

\section{Complete intersections}\label{An application to complete intersections}

In this section we focus on the study of complete intersection varieties of $\P^N$ which lie on some Veronese variety $V_{n,d}$ for $n=1$ and $n=2$; the case $n=3$ and $d=2$ is treated in Proposition \ref{propositionV32}. For generalities on complete intersections we refer to \cite{Commutativealgebra}. To compare with similar  existence and non-existence results for complete intersection on hypersurfaces and their applications we refer to \cite{CarliniChiantiniGeramitaCI,SebastianTripathiUlrich}.

\eatit{\begin{defn}\rm
Let $\X$ be a projective variety in $\P^n$ and consider $\bI(\X)=(f_1,f_2,\dots,f_r)$. We say that $\X$ is a \textit{complete intersection} if $\codim\X=r$ and in this case we write $\X=\mathfrak{C}(f_1,f_2,\dots,f_r)$.
\end{defn}}
\begin{defn}\rm
Let $\X$ be a projective variety, let $\bI(\X)=(f_1,f_2,\dots,f_r)$ be its defining ideal, and let $a_i=\deg f_i$. We say that $\X$ is a \textit{complete intersection of type $(a_1,\ldots,a_r)$} if $\codim\X=r$. 
\end{defn}

\eatit{\begin{prop}\label{unoedue}
If $\X\subseteq V_{n,d}\subseteq\P^N$ is a complete intersection $\X=\mathfrak{C}(f_1,\dots,f_r)$ then there exist $i_1\in\N$ such that $\deg f_{i_1}=1$. Moreover, either $\deg f_i=1$ for each $i=1,\dots,r$ or there exists $i_2\in\N$ such that $\deg f_{i_2}=2$.
\end{prop}
\proof
We set $\deg f_i=a_i$ and we can suppose that $a_1\leq a_2\leq \dots\leq a_r$, so that there exists $p,q\geq 0$ such that $a_1=a_2=\dots=a_p=1$ and $a_{p+1}=a_{p+2}=\dots=a_{p+q}=2$. Since $\X\subseteq V_{n,d}$, we have $\bI(V_{n,d})\subseteq \bI(\X)=(f_1,\dots,f_r)$ so that the generators of $\bI(V_{n,d})$ are generated by $f_1,\dots,f_r$. In particular, since the generators of $V_{n,d}$ are quadrics, we have that $\bI(V_{n,d})\subseteq(f_1,\dots,f_{p+q})$. Suppose by contradiction that $p=0$. Up to a linear change of coordinates, there exists $i_0$ such that $f_1,\dots,f_{i_0}$ are the generators of $\bI(V_{n,d})$. On the one hand, since $V_{n,d}$ is a determinantal variety, $(f_1,\dots,f_{i_0})$ have syzygies $(g_1,\dots,g_{i_0})$ with $\deg g_i=1$ for each $i=1,\dots,i_0$. On the other hand, observing that $(f_1,\dots,f_{i_0})$ is a regular sequence, we get that the only syzygies $(g_1,\dots,g_{i_0})$ of $(f_1,\dots,f_{i_0})$ are the trivial ones. In particular we have $\deg g_i=2$ for each $i=1,\dots,i_0$ so obtaining a contradiction and proving that $p\geq 1$. Finally, we observe that if the generators of $\bI(V_{n,d})$ were generated by $f_1,\dots,f_p$, then $V_{n,d}$ would contain a linear space of $\P^N$ of dimension $N-p$ and this is a contradiction except for $p=N$, because $V_{n,d}$ does not contain linear spaces of dimension $\geq1$. Hence, either $\deg f_i=1$ for each $i=1,\dots, r$ or there is at least a quadric among the $f_1,\dots,f_r$, that is $q\geq 1$. \prfend}

\begin{prop}\label{unoedue}
If $\X\subseteq V_{n,d}\subseteq \P^N$ is a reduced complete intersection of type $(a_1,\dots,a_r)$, with $a_1\leq\ldots\leq a_r$, then $a_1=1$. Moreover, either $r=N$ and $a_1=\ldots=a_N=1$ or  $a_i=2$ for some $i$.
\end{prop}
\proof
Since $\X\subseteq V_{n,d}$, we have that $I(\X)\supseteq I(V_{n,d})$. Recall that $I(V_{n,d})$ is generate by the $2\times 2$ minors of a matrix of linear forms, see \cite{Puccicatalecticant}. Thus, either $a_1=1$ or $a_1=2$. Assume by contradiction that $a_1=2$, that is $a_i>1$ for all $i$, and let $I(\X)=(f_1,\ldots,f_r)$. Thus it is possible to assume that $I(V_{n,d})=(f_1,\ldots,f_q)$ for some $q$. Thus any syzygy of the generators of $I(V_{n,d})$ gives a syzygy of the generators of $I(\X)$. Thus, the determinantal representation of $I(V_{n,d})$ yields that the generators of $I(\X)$ have a syzygy of linear forms. Hence a contradiction since $f_1,\ldots,f_r$ are a regular sequence and their syzygies, given by the Koszul complex, only contain elements of degree at least $a_1=2$. In conclusion, $a_1=1$.
To complete the proof, we let $p$ be the largest index such that $a_p=1$. Assume by contradiction that $p< r$ and $a_{p+1}>2$. Thus, $I(V_{n,d})\subseteq (f_1,\ldots,f_p)$ and the latter is the ideal of a linear space $\Lambda$ of codimension $p$ such that $\Lambda\subseteq V_{n,d}$. Hence a contradiction since no Veronese variety contains a positive dimensional linear spaces. In conclusion, either $p=r=N$ or $p<r$ and $a_{p+1}=2$.
\prfend
\noindent As a straightforward application of this proposition we get the following result.
\begin{thm}\label{cicurve}
Let $\X\subseteq V_{1,d}\subseteq\P^d$ be a reduced complete intersection with $d\geq 3$. Then either $\X$ is a point or $\X$ is a set of two points.
\end{thm}
\proof
Since $V_{1,d}$ is not a complete intersection and $\dim V_{1,d}=1$, it must be $\dim \X=0$. Since $\X$ is a complete intersection, it follows from Proposition \ref{unoedue} that $\X$ is contained in the intersection of a hyperplane and $V_{1,d}$. In particular we have $|\X|\leq d$ and using Remark \ref{hfcurve} we get the Hilbert function of $\X$:
$$
\begin{array}{c|cccc}
t & 0 & 1 & 2\ & \dots\\
\hline
H_\X(t) & 1 & |\X| & |\X| & \dots\\
\Delta H_\X(t) & 1 & |\X|-1 & 0 & \dots
\end{array}
$$
Finally, since $\X$ is a complete intersection, $\Delta H_\X (t)$ is the Hilbert function of an artinian Gorenstein ideal. Thus it is symmetric and hence either $\Delta H_\X(1)=0$ or $\Delta H_\X(1)=1$, that is $|\X|=1$ or $|\X|=2$. \prfend

\begin{remark}\rm
We note that the proof of Theorem \ref{cicurve} also shows that any Gorenstein reduced zero-dimensional subscheme of $V_{1,d}$ must be degenerate.
\end{remark}
We can now finally describe reduced complete intersection subvarieties of Veronese surfaces.
\begin{thm}\label{cisuperfici}
If $\X\subseteq V_{2,d}\subseteq \P^N$ is a reduced complete intersection of type $(a_1,\dots,a_r)$, with \linebreak $a_1\leq\dots\leq a_r$ then one of the following holds:
\begin{enumerate}
	\item $(d,r,(a_1,a_2,\dots,a_r))=(2,4,(1,1,1,2))$, that is $\X$ is a conic lying on $V_{2,2}$;
	\item $(d,r,(a_1,a_2,\dots,a_r))=(2,5,(1,1,1,2,a_5))$, any $a_5\in\N$, that is $\X$ is a set of $2a_5$ complete intersection points of a conic lying on $V_{2,2}$ and a hypersurface of degree $a_5$;
	\item\label{case3thmcisurfaces} $(d,r,(a_1,a_2,\dots,a_r))=(d,N,(1,1,\dots,1))$ for any $d\geq 2$, that is $\X$ is a reduced point;
	\item $(d,r,(a_1,a_2,\dots,a_r))=(d,N,(1,1,\dots,1,2)$ for any $ d\geq 2$, that is $\X$ is a set of two reduced points.
\end{enumerate}
\end{thm}
\proof
Let $\bI(\X)=(f_1,\dots,f_r)$ be the ideal of $\X$. We can suppose without loss of generality that $a_1\leq a_2\leq\dots\leq a_r$. Moreover, we let $p$ be the number of $a_i$ equal to 1 and $q$ be the  number of $a_i$ equal to 2, that is $a_1=a_2=\dots=a_p=1$, $a_{p+1}=a_{p+2}=\dots=a_{p+q}=2$. By Proposition \ref{unoedue} we have that $p\geq 1$. If $p=N$ we are trivially in case \eqref{case3thmcisurfaces} thus from now on we suppose $p<N$ and thus, using again Proposition \ref{unoedue}, we also have $q\geq 1$.
Since $\X\subseteq V_{2,d}$, $\Delta H_{\X}(t)$ must satisfy conditions \eqref{conduno} and \eqref{conddue} of Theorem \ref{main}. We now use condition \eqref{conduno} to rule out some cases. We start by computing $\Delta H_{\X}(1)$ and $\Delta H_{\X}(2)$. We have that
$$H_\X(1)=\dim\K[y_0,\dots,y_N]-\dim (f_1,\dots,f_r)_1=N-p+1$$
$$\Delta H_{\X}(1)=H_{\X}(1)-H_\X(0)=N-p.$$
Now, let us consider the map
$$\begin{array}{ccc}
 \K[y_0,\dots,y_N] & \to & \K[y_0,\dots,y_N]/(f_1,\dots,f_p)\cong\K[y_0,\dots,y_{N-p}]\\[1ex]
  f & \mapsto & \tilde{f}
\end{array}.$$
Since $(f_1,f_2,\dots,f_r)$ is a regular sequence we have that
$$H_\X(2)=\dim\big(  \K[y_0,\dots,y_N]/(f_1,\dots,f_r)\big)  _2=\dim\big( \K[y_0,\dots,y_{N-p}]/(\tilde{f}_{p+1},\dots,\tilde{f}_{r})\big) _2$$
$$=\dim{\K[y_0,\dots,y_{N-p}]_2}-q={N-p+2\choose 2}-q$$
and hence
$$\Delta H_\X(2)={N-p+2\choose 2}-q-(N-p+1)={(N-p+1)(N-p)\over 2}-q.$$
Note that $H_\X(1)=N-p+1<N+1=H_{V_{2,d}}(1)$, thus $t_1=0$. Since $t_1=0$ we have that
$$\mu_2(d,0,N-p)=\floor[\bigg]{{2d+3-\sqrt{1+8p}\over 2}},$$
so that by \eqref{conduno} we get
\begin{equation}\label{cattiva}
\floor[\bigg]{{2d+3-\sqrt{1+8p}\over 2}}\geq \ceil[\bigg]{{(N-p+1)(N-p)\over 2d}-{q\over d}}.
\end{equation}
Since $q\leq N-p$, if $(p,d,q)$ is a solution of \eqref{cattiva}, then $(p,d)$ is a solution of
$$
{2d+3-\sqrt{1+8p}\over 2}\geq {(N-p)(N-p-1)\over 2d}
$$
and setting $\alpha=N-p$ the previous inequality yields
\begin{equation}\label{cattiva2}
\underbrace{2d+3-\sqrt{(2d+3)^2-8(\alpha+1)}}_{f(\alpha)}\geq \underbrace{{\alpha(\alpha-1)\over d}}_{g(\alpha)}.
\end{equation}
A standard calculus argument on $f(\alpha)$ and $g(\alpha)$ shows  that if $(\alpha,d)$ is a solution of \eqref{cattiva2} then $\alpha\leq 3$. As a consequence, if $(p,d,q)$ is a solution of \eqref{cattiva} then $p\geq N-3$, thus one has to look for solutions $(p,d,q)$ only for $p=N-3,N-2,N-1$. Since $r\leq N$, one can easily check that the solutions of (\ref{cattiva}) are the following:
\begin{itemize}[leftmargin=*]
\item $(p,d,q)=(N-1,d,1)$ for each $d\geq 2$.\\
In this case $\X$ is a complete intersection of type $(a_1,\dots,a_{N})$ with $a_1=\dots=a_{N-1}=1$ and $a_N=2$, that is $\X$ is a set of two reduced points and we are in case \textit{4.} 
\item $(p,d,q)=(3,2,1)$.\\
In this case $N=5$ and we know that $a_1=a_2=a_3=1$ and $a_4=2$, thus we distinguish two subcases:
\begin{itemize}[leftmargin=*]
	\item If $r=4$ then $\X$ is a complete intersection of type $(1,1,1,2)$, that is $\X$ is a conic lying on $V_{2,2}$ and we are in case \textit{1.}
	\item If $r=5$ then $\X$ is a complete intersection of type $(1,1,1,2,a_5)$, that is $\X$ is the intersection of a conic lying on $V_{2,2}$ with a hypersurface of degree $a_5$ and we are in case \textit{2.}
\end{itemize}
\item $(p,d,q)=(3,2,2)$.\\
In this case $\X$ is a complete intersection of type $(1,1,1,2,2)$, thus it this is a special case of \textit{2.}
\item $(p,d,q)=(6,3,3)$.\\
In this case $\X$ is a complete intersection of type $(1,1,1,1,1,1,2,2,2)$. We want to show that such $\X$ can not lie on $V_{2,3}$. Let us suppose by contradiction that $\X\subseteq V_{2,3}$ and set $\Y=\nu_{2,3}^{-1}(\X)$. From Lemma \ref{Imphi} it follows that
$$H_\Y(3)=H_\X(1)=N-p+1=9-6+1=4.$$
Thus, the only way to complete the gaps of $H_\Y$ is
$$\begin{array}{c|ccccccccccc}
t & 0 & 1 & 2 & 3 & 4 & 5 & 6 & 7 & 8 &\dots\\ \hline
H_\Y(t) & 1 & 2 & 3 & 4 & 5 & 6 & 7 & 8 & 8 & \dots
\end{array}$$
and this shows that $\Y$ lies on a line. As a consequence, $\X$ lies on a rational normal curve. Hence, since $\X$ is a complete intersection and $|\X|=8$, this is a contradiction by Theorem \ref{cicurve}.
\end{itemize}
The result is now proved.
\prfend

\section{More results and open problems}\label{More results and open problems}

In this section, we show some more results about complete intersections on Veronese varieties. A complete characterization of the Hilbert functions of subvarieties of Veronese varieties of dimension larger than two seems, presently, to be out of reach. Nevertheless, we can deal with the case of the threefold $V_{3,2}$. This case leads us to formulate Conjecture \ref{theconjecture}.

\begin{prop}\label{propositionV32}
Let $\X\subseteq V_{3,2}\subseteq\P^9$ be a reduced subvariety. Then $\X$ is a complete intersection of type $(a_1,\ldots,a_r)$, with $a_1\leq\dots\leq a_r$ if and only if $\X$ is one of the following
\begin{itemize}
    \item $r=9$,$a_1=\ldots=a_9=1$, that is $\X$ is a reduced point;
    \item $r=9$,$a_1=\ldots=a_8=1,a_9=2$, that is $\X$ is a set of two reduced points;
    \item $r=9$,$a_1=\ldots=a_7=1,a_8=2,a_9=b$, any $b\geq 2$, that is $\X=\mathcal{C}\cap H_b$ for $\mathcal{C}\subseteq V_{3,2}$ a conic and $H_b$ a degree $b$ hypersurface;
    \item $r=8$,$a_1=\ldots=a_7=1,a_8=2$, that is $\X$ is a conic.
\end{itemize}
\end{prop}
\proof
First we consider the case $\dim\X=0$. Let $\X\subseteq V_{3,2}$ be a reduced complete intersection of type $(a_1,\dots,a_9)$ with $a_1\leq a_2\leq \dots\leq a_9$. Also let $p,q\in\N$ such that $a_1=\dots=a_p=1$ and $a_{p+1}=\dots=a_{p+q+1}=2$. By Proposition \ref{unoedue} either $p=9$ or $p\geq 1$ and $q\geq 1$. If $p=9$ then $\X$ is just a reduced point, thus from now on we suppose $p,q\geq 1$. By the same argument used in the proof of Theorem \ref{cisuperfici} we get
$$H_{\X}(1)=10-p, \quad H_{X}(2)={11-p\choose 2}-q={p^2-21p+110\over 2}-q.$$
Since $\X\subseteq V_{3,2}$, by Lemma \ref{Imphi} there exists $\Y=\nu_{3,2}^{-1}(\X)\subseteq\P^3$ such that  $H_{\X}(t)=H_{\Y}(td)$ for all $t\geq 0$. In particular, we have that
$$H_{\Y}(2)=10-p, \quad H_{\Y}(4)={p^2-21p+110\over 2}-q.$$
Now fix $1\leq p\leq 8$. By Macaulay's theorem for Hilbert functions (see \cite{StanleyHF}, Theorem 2.2) it follows that if the 3-binomial expansion of $H_{\Y}(2)$ is
$$H_{\Y}(2)={m_2\choose 2}+{m_1\choose 1}$$
where $m_2>m_1$, then 
$$H_{\Y}(4)\leq {m_2+2\choose 4}+{m_1+2\choose 3}=M(p)$$
On the other hand, since $1\leq q\leq 9-p $, we have that 
$$H_{\Y}(4)\geq {p^2-21p+110\over 2}-(9-p)=m(p). $$
Thus, if $M(p)-m(p)<0$, then $\X$ does not exist. Computing we get the following table
\begin{table}[H]
\centering
\begin{tabular}{c|c|c|c|c}
$p$ & $H_{\Y}(2)$ & $M(p)$ & $m(p)$ & $M(p)-m(p)$ \\ \hline
1   & 9           & 25     & 37     & -12          \\
2   & 8           & 19     & 29     & -10          \\
3   & 7           & 16     & 22     & -6           \\
4   & 6           & 15     & 16     & -1           \\
5   & 5           & 9      & 11     & -2           \\
6   & 4           & 6      & 7      & -1           \\
7   & 3           & 5      & 4      & 1          \\
8   & 2           & 2      & 2      & 0          
\end{tabular}
\end{table}
\noindent Hence, $\X$ is either of type $(1,1,1,1,1,1,1,2,a_9)$ or of type $(1,1,1,1,1,1,1,1,2)$. In the first case $\X$ is a set of $2a_9$ reduced points lying on a conic $\bC\subseteq V_{3,2}$ and in the second case $\X$ is a set of 2 reduced points. The discussion for reduced 0-dimensional complete intersection is now completed.
Now let $\X\subseteq V_{3,2}$ be a positive dimensional reduced complete intersection  of type $(a_1,\dots,a_r)$. For any choice of integers $a_{r+1},\dots, a_9$,  we can choose suitable hypersurfaces $H_{a_i}$ of degree $a_i$ in such a way that 
$$
\X'=\X\cap H_{a_{r+1}}\cap\ldots\cap H_{a_9}
$$
is a complete intersection of type $(a_1,\ldots,a_r,a_{r+1},\ldots,a_9)$.
Moreover, we can choose the degrees $a_i$ in a such a way that
\[
a_i\leq a_{i+1}  
\]
for all $i$ and $3\leq a_{r+1}$.
Thus $\X'\subseteq V_{3,2}$ is a zero dimensional complete intersection of type $(a_1,\ldots,a_9)$. As a consequence, since we can freely choose the degrees $a_i$ for $i\geq r+1$ we have that $r=8$ and 
\[
(a_1,\dots,a_8,a_9)=(1,1,1,1,1,1,1,2,a_9).
\]
 Hence, $\X$ is of type $(1,1,1,1,1,1,2)$, that is $\X$ is a conic.
\prfend
\noindent We can now state the following conjecture, which we already proved for $n\leq 2$ any $d$ and for $n=3$ and $d=2$, see Proposition \ref{propositionV32} and Theorems \ref{cicurve} and \ref{cisuperfici}.

\begin{conj}\label{theconjecture}
Let $\X\subseteq V_{n,d}\subseteq\P^N$ be a reduced subvariety with $d>1$. Then $\X$ is a complete intersection of type $(a_1,\ldots,a_r)$, with $a_1\leq\dots\leq a_r$ if and only if
\begin{itemize}
    \item $r=N$,$a_1=\ldots=a_N=1$, any $n,d$, that is $\X$ is a reduced point;
    \item $r=N$,$a_1=\ldots=a_{N-1}=1$,$a_N=2$, any $n,d$, that is $\X$ is a set of two reduced points;
    \item $r=N$,$a_1=\ldots=a_{N-2}=1$,$a_{N-1}=2$,$a_N=b$, any $n$, $d=2$, any $a\geq 2$, that is $\X=\mathcal{C}\cap H_b$ for $\mathcal{C}\subseteq V_{n,2}$ a conic and $H_b$ a degree $b$ hypersurface;
    \item $r=N-1$,$a_1=\ldots=a_{N-2}=1,a_{N-1}=2$, $d=2$, any $n$, that is $\X$ is a conic.
\end{itemize}
\end{conj}
In the case of the Veronese threefold $V_{3,2}$, see proof of Proposition \ref{propositionV32}, the complete knowledge of the zero dimensional complete case allows us to complete the proof. This is true in general as shown by the following Lemma.
\begin{lem}
If Conjecture \ref{theconjecture} holds for all reduced zero dimensional subvariety of $V_{n,d}$, then it holds for all reduced subvarieties of $V_{n,d}$.
\end{lem}
\proof
Let $\X\subseteq V_{n,d}\subseteq\P^N$ be a reduced complete intersection of type $(a_1,\ldots,a_r)$ with $r<N$, that is $\X$ is positive dimensional. Then, for any choice of integers $ a_{r+1}, \ldots, a_N$ we can choose suitable hypersurfaces $H_{a_i}$ of degree $a_i$ in such a way that
\[
\X'=\X\cap H_{a_{r+1}}\cap\ldots\cap H_{a_N}
\]
is a complete intersection of type $(a_1,\ldots,a_r,a_{r+1},\ldots,a_N)$. Moreover, we can choose the degrees $a_i$ in a such a way that
\[
a_i\leq a_{i+1}  
\]
for all $i$ and $3\leq a_{r+1}$.
Thus $\X'\subseteq V_{n,d}$ is a zero dimensional complete intersection of type $(a_1,\ldots,a_N)$. Since we are assuming that the conjecture holds for such an $\X'$ and since we can freely choose the degrees $a_i$ for $i\geq r+1$ we have that $r=N-1$ and
\[
(a_1,\ldots,a_{N-1},a_N)=(1,\ldots,1,2,a_N)
\]
and thus $d=2$ and $\X$ is a conic.
Hence the conjecture holds for $\X$.
\prfend
\appendix
\section{Appendix}
We now prove  Lemma \ref{tec}.
\dd 
\proof We distinguish four cases depending on the value of $p$. For each of them we give a function $\tilde{h}(i)$ satisfying \textit{1.} and \textit{2.} and such that 
$$\tilde{h}(d)=\begin{cases}
\floor[\big]{{2d(t+1)+3-\sqrt{1+8p}\over 2}}, &\text{if } 1\leq p\leq {d+1\choose 2}\\
dt-n, & \text{if } {d+1\choose 2}+dn<p\leq {d+1\choose 2}+d(n+1), 0\leq n\leq dt
\end{cases}$$
Then we show that for any function $h'(i)$ satisfying \textit{1.} and \textit{2.} it holds that $h'(d)\leq \tilde{h}(d)$. We do this in detail in case \eqref{appendice1} and for the remaining cases we produce the function $\tilde{h}$.
\begin{enumerate}[leftmargin=*]\label{appendice1}
	\item $p={n\choose 2},\; 1\leq n\leq d+1$
	\dd
	In this case we set
	$$\tilde{h}(i)=\begin{cases}
	dt+i+1, & \text{if } 1\leq i\leq d-n \\
	d(t+1)-n+2, & \text{if }  d-n+1\leq i\leq d
	\end{cases}.$$
	We have 
	$$\sum_{i=1}^{d}\tilde{h}(i)=d^2t+{d(d+3)\over 2}-{n(n-1)\over 2}=d^2t+{d(d+3)\over 2}-p=\sum_{i=1}^{d}h(i)$$
	and
    $$\tilde{h}(d)=d(t+1)-n+2=\floor[\bigg]{{2d(t+1)+3-\sqrt{1+8p}\over 2}}$$	
	hence $\tilde{h}(i)$ is as we want.
	Now let us suppose that there exists $h'(i):\left\lbrace 1,2,\dots,d\right\rbrace \to\N$ satisfying \textit{1.} and \textit{2.} and such that $h'(d)>\tilde{h}(d)$, that is $h'(d)=d(t+1)-n+2+a,\, a\geq 1$. Since $n\leq d+1$, we have $h'(d)\geq dt+1+a\geq dt+2$. As a consequence (observe that $dt+2$ is the maximum value of $h(1)$), by \textit{1.} and \textit{2.}, it follows that $h'(i)$ is increasing at least until reaching the value $h'(d)$. In particular, if we set $i'=\min\left\lbrace 1\leq i\leq d\;|\; h(i)=d(t+1)-n+2+a\right\rbrace $ we have
	$$dt+i'+1=d(t+1)-n+2+a$$
	so that $i'=d-n+1+a$ and $i_0\geq i'$. Hence, using again \textit{1.} and \textit{2.}, we get
	$$
	\begin{array}{cc}h'(i)=\tilde{h}(i)=dt+i+1 &\text{if } i\leq d-n\\[2ex]
	h'(i)>\tilde{h}(i) & \text{if } i\geq d-n+1
	\end{array}$$
Hence:
$$\sum_{i=1}^d h'(i)=\sum_{i=1}^{d-n}\tilde{h}(i)+\sum_{i=d-n+1}^{d}\underbrace{h'(i)}_{>\tilde{h}(i)}>\sum_{i=1}^d \tilde{h}(i)=\sum_{i=1}^d h(i)$$
and this is a contradiction.
\item ${n\choose 2}<p<{n+1\choose 2},\; 1\leq n\leq d$
\dd
Let $b\in\Z$ be such that $p={n\choose 2}+b$. In this case we set
$$\tilde{h}(i)=\begin{cases}
dt+i+1, & \text{if } 1\leq i\leq d-n+1\\
d(t+1)-n+2, & \text{if } d-n+2\leq i\leq d-b\\
d(t+1)-n+1, & \text{if } d-b+1\leq i\leq d
\end{cases}.$$

\item $p={d+1\choose 2}+d(n+1),\, 0\leq n\leq dt$\\
In this case we set 
$$\tilde{h}(i)=dt-n,\; 1\leq i\leq n.$$
We have
$$\sum_{i=1}^d\tilde{h}(i)=d(dt-n)=d^2t-dn=d^2t+{d(d+3)\over 2}-p=\sum_{i=1}^dh(i)$$
and 
$$\tilde{h}(d)=dt-n$$
hence $\tilde{h}(i)$ is as we want. 
\item ${d+1\choose 2}+nd<p<{d+1\choose 2}+(n+1)d,\, 1\leq n\leq dt$\\
Let $b\in\Z$ be such that $p={d+1\choose 2}+nd+b$. In this case we set
$$\tilde{h}(i)=\begin{cases}
dt+1-n, & \text{if } 1\leq i\leq d-b\\
dt-n, & \text{if } d-b+1\leq i\leq d
\end{cases}$$
\end{enumerate}\prfend
\eatit{\begin{lem}\label{riseq}\rm
The integer solutions of the equation \eqref{cattiva} are
\begin{itemize}
	\item $(p,d)=(N,d)$ for each $d\geq 2$;
	\item $(p,d)=(N-1,d)$ for each $d\geq 2$;
	\item $(p,d)=(3,2)$;
	\item $(p,d)=(6,3)$.
\end{itemize}
\end{lem}
\proof
First of all we observe that since $q\leq N-p$ therefore if $(p,d)$ is a solution of \eqref{cattiva} then it is also a solution of
$$
{2d+3-\sqrt{1+8p}\over 2}\geq {(N-p)(N-p-1)\over 2d}.
$$
By setting $\alpha:=N-p$ the previous inequality can be rewritten as follows
\begin{equation}
\underbrace{2d+3-\sqrt{(2d+3)^2-8(\alpha+1)}}_{f(\alpha)}\geq \underbrace{{\alpha(\alpha-1)\over d}}_{g(\alpha)}.
\end{equation}
Since $f'(\alpha)\approx {1\over\sqrt{\alpha}}$ and $g'(\alpha)\approx\alpha$, if there exists $\alpha_0$ such that $f'(\alpha_0)>g'(\alpha_0)$ then $g'(\alpha)>f'(\alpha)$ for each $\alpha\geq\alpha_0$. As a consequence, if there exists $\alpha_1\geq a_0$ such that $g(\alpha_1)>f(\alpha_1)$ then $g(\alpha)>f(\alpha)$ for each $\alpha\geq \alpha_1$. A straightforward computation shows that for each $d\geq2$ it holds that $g(4)>f(4)$ and $g'(4)>f'(4)$ and hence if $(\alpha,d)$ is a solution of \eqref{cattiva} then $\alpha\leq 3$. It is easy to check that $(0,d)$ and $(1,d)$ are solutions of \eqref{cattiva} for each $d\geq 2$ while $(2,d)$ is a solution if and only if $d=2$ and $(3,d)$ is a solution if and only if $d=3$.\prfend}

\newcommand{\etalchar}[1]{$^{#1}$}

\end{document}